\documentclass[12pt]{amsart}
\usepackage{amsmath}

\def\ol{\overline}
\def\e{\epsilon}
\def\lf{\left}
\def\ri{\right}

\def\wt{\widetilde}

\def\p{\partial}

\newcommand\ce{{\mathbb C}}
\newcommand\C{{\mathbb C}}
\def\ibar{{\bar\imath}}

\def\jbar{{\bar\jmath}}

\def\K{K\"ahler }
\def\KR{K\"ahler-Ricci }
\def\KRF{K\"ahler-Ricci flow }
\def\KRS{K\"ahler-Ricci soliton }
\def\KRS{K\"ahler-Ricci soliton }

\def\be{\begin{equation}}
\def\ee{\end{equation}}
\def\ol{\overline}
\def\lf{\left}
\def\ri{\right}

\def\e{\epsilon}

\def\ijb{{i\jbar}}

\def\wt{\widetilde}

\def\p{\partial}

\def\C{\Bbb C}
\def\P{\Bbb P}
\def\cn{\Bbb C^n}
\def\wh{\widehat}

\def\wt{\widetilde}

\def\p{\partial}

\def\C{\Bbb C}
 
\def\KRF{K\"ahler-Ricci flow }

\newtheorem{thm}{Theorem}[section]

\newtheorem{lem}{Lemma}[section]
\newtheorem{prop}{Proposition}[section]
\newtheorem{cor}{Corollary}[section]
\theoremstyle{definition}

\theoremstyle{remark}
\newtheorem{rem}{Remark}[section]
\numberwithin{equation}{section}
 \begin{document}
\title{Non-negatively curved K\"ahler manifolds with average quadratic curvature decay}

\author{Albert Chau}
\address{Waterloo University, Department of Pure Mathematics,
  200 University avenue, Waterloo, ON N2L 3G1, CANADA}
\email{a3chau@math.uwaterloo.ca}

\author{Luen-Fai Tam$^1$}

\thanks{$^1$Research
partially supported by Earmarked Grant of Hong Kong \#CUHK403005}

\address{Department of Mathematics, The Chinese University of Hong Kong,
Shatin, Hong Kong, China.} \email{lftam@math.cuhk.edu.hk}

\renewcommand{\subjclassname}{%
  \textup{2000} Mathematics Subject Classification}

\date{Sept.  2005}



\begin{abstract}
Let $(M, g)$ be a   complete non compact \K manifold with
non-negative and bounded holomorphic bisectional curvature.
Extending our techniques developed in \cite{CT3}, we prove that
the universal cover $\wt M$ of  $M$ is biholomorphic to $\ce^n$
provided either that $(M, g)$ has average quadratic curvature
decay, or $M$ supports an eternal solution to the \KR flow with
non-negative and uniformly bounded holomorphic bisectional
curvature. We also classify certain local limits arising from the
\KR flow in the absence of uniform estimates on the injectivity
radius.
\end{abstract}

\maketitle \markboth{Albert Chau and Luen-Fai Tam} {Non-negatively curved \K manifolds with average quadratic curvature decay}

\section{Introduction}

The relation between the geometry and the complex structure of a
complete \K manifold with positive holomorphic bisectional
curvature has long been studied.  It has been conjectured that
such a simply connected \K manifold $M^n$ is biholomorphic to
either $\ce \P^n$ or $\ce^n$ depending on whether $M^n$ is
respectively compact or non-compact. If $M^n$ is compact the
conjecture was due to Frankel and was later independently solved by Mori \cite{M} and
Siu-Yau \cite{SY}.  If $M^n$ is non-compact, the conjecture
is due to Yau \cite{Y} and remains unsolved.  Several authors have worked on Yau's conjecture and we refer to \cite{CT3} for a review of past results.   In
\cite{CT3} the authors proved the following result supporting
Yau's conjecture

\begin{thm}\label{s1t0}
Let $(M^n,\wt g)$ be a complete noncompact K\"ahler manifold
with nonnegative and bounded holomorphic bisectional
curvature and maximal volume growth.  Then $M$ is biholomorphic to $\C^n$.
\end{thm}
Here maximum volume growth means that

\begin{equation}\label{maxvol} Vol(B(p, r)) \geq C_1r^{2n};\hspace{11pt} \forall r\in [0,
  \infty)\end{equation}  for some  $C_1>0$ and  $p \in M$.  Maximum volume growth of such a \K manifold is closely related to the
 the average quadratic curvature decay condition
\begin{equation}\label{quadraticdecay}\frac{1}{V_x(r)}\int_{B_x(r)}R\le \frac{C_2}{1+r^2}\end{equation}
for some $C_2>0$, all $x\in M$ and all $r>0$. Here $B_x(r)$ is the
geodesic ball around $x$ with volume $V_x(r)$ and $R$ is the
scalar curvature of $M$.  It was conjectured by Yau that
(\ref{maxvol}) implies (\ref{quadraticdecay}) for a complete
non-compact \K manifold with non-negative holomorphic bisectional
curvature.  Provided that the curvature is bounded, this was
recently confirmed by Ni \cite{Ni2} (this was earlier confirmed by
Chen-Tang-Zhu \cite{CTZ}  for the case of dimension 2 and Chen-Zhu
\cite{CZ3} in all dimensions under the additional condition that
the curvature operator is nonnegative). In general, (1.2) does not
imply (1.1).  However, one may expect that Theorem \ref{s1t0}
  is still true if we replace the maximal volume growth
condition with average quadratic curvature decay.  We confirm this
expectation in the following

\begin{thm}\label{s1t1}
Suppose $(M^n, g)$ has holomorphic bisectional curvature which is
bounded, non-negative and has average quadratic curvature decay.
 Then $M$ is holomorphically covered by
$\ce^n$.
\end{thm}
Theorem \ref{s1t1} was proved by the authors in \cite{CT3} under the additional
assumption that the curvature operator is non-negative.

As in \cite{CT3}, we will use the \KR flow:
\begin{equation}\label{s1e0}
\begin{cases}
    \frac{\p }{\p t}\tilde {g}_{i\jbar}(x,t)& =- \tilde{R}_{i\jbar}(x,t);  \\
    \tilde{g}_{i\jbar}(x,0)&  = {g}_{i\jbar}(x).
\end{cases}
\end{equation}

For $(M, g)$ as in Theorem \ref{s1t1}, it is now well known by \cite{Sh2} (see
also \cite{NT}) that (\ref{s1e0}) has a long time solution $\tilde g(t)$, $0\le
t<\infty$.  If we let $g(x,t)=e^{-t}\tilde g(x,e^t)$, then we obtain a solution
to the normalized \KR flow
\begin{equation}\label{s1e1}
\frac{\p }{\p t}g_{i\jbar}(x,t)=-  R_{i\jbar}(x,t)-g_\ijb(x,t)
\end{equation}
 for $-\infty<t<\infty$.  Moreover $g(t)$ has uniformly bounded nonnegative holomorphic bisectional curvature on
$M\times(-\infty,\infty)$ by \cite{Sh0}, \cite{Sh2} and \cite{Ha4}, and for any
$k\ge 1$ the norms $|D^kRm|$ of the covariant derivatives of the curvature tensor
of $g(t)$ with respect to $g(t)$ are uniformly bounded on $M\times [a,\infty)$
for any $a>-\infty$.

We will also prove a uniformization theorem for {\it eternal} solutions to the
\KR flow.  Recall that an    eternal  solution of the \KR flow on a complex
manifold $M$ is a smooth family of complete \K metrics $g(t)$ on $M$ satisfying
\begin{equation}\label{s1e2}
\frac{\p }{\p t} {g}_{i\jbar}(x,t)=- {R}_{i\jbar}(x,t)
\end{equation}
for all $t\in (-\infty,\infty)$.  We prove that
\begin{thm}\label{s1t2}
Let $(M, g(t))$ be a complete eternal solution to (\ref{s1e2})
such that for all $t$, $g(t)$ has non-negative holomorphic
bisectional curvature which is uniformly bounded on $M$
independent of $t$.  Then $M$ is holomorphically covered by
$\ce^n$.
\end{thm}
As before, by \cite{Sh0} and \cite{Ha4}, $|D^kRm|$ is also
uniformly bound on $M\times[a,\infty)$ for $g(t)$ for any
$a>-\infty$.

\begin{rem} By the results of \cite{MSY} (see also
\cite{Chen,NT2}), if $M$ is complete noncompact with bounded nonnegative
bisectional curvature and if the curvature decays faster than quadratic in the
average sense, then $M$ is flat.  Hence Theorem 1.2 addresses the maximal
(quadratic) curvature decay case for nonflat $M$. On the other hand,  by the
Harnack inequality \cite{cao} and the decay estimates in \cite[section 6]{Sh1}
(see also \cite[Corollary 2.1]{NT}, it is seen that the average curvature of $(M,
g)$ in Theorem 1.3 cannot decay faster than linearly uniformly at all points.  In
this sense, Theorem 1.3 addresses the case of minimal (linear) curvature decay by
the results in \cite{NT2}.
\end{rem}

By comparing (\ref{s1e1}) and (\ref{s1e2}), we may combine Theorems \ref{s1t1}
and \ref{s1t2} in the following:
\begin{thm}\label{s1t3} Let $M^n$ be a noncompact complex manifold.
Suppose there is a smooth family of complete \K metrics $g(t)$ on $M$ such that
for $\kappa=0$ or $1$, $g(t)$ satsifies
\begin{equation}\label{s1e3}
\frac{\p }{\p t} {g}_{i\jbar}(x,t)=- {R}_{i\jbar}(x,t)-\kappa
g_\ijb(x,t)
\end{equation}
for all $t\in (-\infty,\infty)$ such that for every $t$, $g(t)$ has uniformly
bounded non-negative holomorphic bisectional curvature on $M$ independent of $t$.
Then $M$ is holomorphically covered by $\ce^n$.
\end{thm}

By the results in \cite{Cao}, if in Theorem \ref{s1t3} we assume that the Ricci
curvature is positive and the scalar curvature attains its maximum in spacetime,
then $(M,g(t))$ is a gradient K\"ahler-Ricci soliton of steady type if
$\kappa=0$, and of expanding type if $\kappa=1$. By the results in \cite{b, CT2},
$M$ is biholomorphic to $\cn$. Hence Theorem \ref{s1t3} can also be considered as
a generalization of the results in \cite{b, CT2}.

On the other hand, it was proved by  Fangyang Zheng, see
\cite[Theorem 5.2]{NT2},   that a complete noncompact \K manifold
with nonnegative sectional curvature and with positive Ricci
curvature at some point is simply connected.    Hence we have the
following:

\begin{cor}\label{s1c1} Let $(M^n,g)$ be as in Theorem \ref{s1t1}.
If addition, $M$ has nonnegative sectional curvature and the Ricci curvature is
positive at some point, then $M$ is biholomorphic to $\cn$.
\end{cor}

\begin{rem}\label{reduction} If we take $\pi:\wh M\to M$ to be the universal holomorphic
covering of $M$ and let $\wh g(t)=\pi^*(g(t))$, then $(\wh M,\wh
g)$ still satisfy the conditions of Theorem \ref{s1t3}. To prove
the theorem, it is sufficient to prove that $\wh M$ is
biholomorphic to $\cn$. By \cite{cao1}, we may further assume that
the Ricci curvature of $\wh g(x,t) $ is positive for all $x$ and
$t$. Hence from now on we assume $M$ in Theorem \ref{s1t3} is
simply connected and $g(t)$ has positive Ricci curvature for all
$t$.
\end{rem}

The Authors would like to thank  Huai-Dong Cao, David Glickenstein, Richard
Hamilton, Lei Ni, Hung-Hsi Wu and Shing-Tung Yau for helpful discussions.

\section{Local limit solution}
Throughout the section, we assume $(M, g(t))$ satisfies the conditions of Theorem
\ref{s1t3}, and that $M$ is simply connected and $g(t)$ has strictly positive
Ricci curvature for all $t$. Fix some point $p\in M$, some time sequence $t_k \to
\infty$ and consider the sequence $(M, g(t_k +t), p)$ of long time solutions to
the \KRF centered at $p$. Suppose the injectivity radius of $g(t)$ at $p$ has a
uniform lower bound. Then by Hamilton's compactness \cite{Ha2}, this sequence has
a convergent subsequence converging to a solution $h(t)$ to the \KRF on a limit
complex manifold $N$. Furthermore, by Cao's classification of limits for \KR flow
\cite{Cao}, this limit must either be a steady or expanding gradient \KR soliton,
depending on whether $\lambda=0$ or $\lambda=1$. In this section, we show that in
the absence of an injectivity radius estimate, we may still have such a soliton
limit, but in a local sense. We will consider a certain locally lifted
subsequence limit of $(M, g(t_k +t), p)$ around $p$. Our first goal will be to
show that this local limit is also either an expanding or steady gradient \KR
soliton in a certain sense (Theorem \ref{s2t1}).  We will then relate this to the
local asymptotic behavior of $g(t)$ at $p$ in our main Theorem \ref{s2t3}. In the
absence of injectivity radius estimates, Glickenstein \cite{Gl} constructed a
global limit solution from a solution to the Ricci flow as above, allowing for
the possibility of Gromov Hausdorff convergence to a limiting metric space of
dimension lower than that of $M$. We refer the reader to \cite{Gl} for details on
the construction of this limit and its application, and in particular to
\cite{CGP1, CGP2} for applications in three dimensions. Our local limit is just
the first step of Glickenstein's construction and in fact depends only on
Proposition \ref{s2p1} and the simple fact that a lifting of a solution to the
Ricci flow is still a solution to the flow.  For recent work relating this and in general, on the existence and classification of
limits to the Ricci flow, we refer to the works of Ye \cite{Ye} and Lott \cite{L}.

Recall that from the time independent bounds on the curvature of
$g(t)$, we have corresponding bounds on all covariant derivatives
of the curvature by the \KR flow. Hence for $t\ge a$ with $a>-\infty$, we may
assume that these bounds on all covariant derivatives of the
curvature of $g(t)$ are also time independent. The proof of
Proposition 1.2 in \cite{TY} then gives (see also
\cite[Proposition 2.1]{CT3}):

\begin{prop}\label{s2p1}
There exist positive constants $r$  and $C$ such that for each
$t\ge-1$ there is a holomorphic map $\Phi_t$ from the Euclidean
  ball $D(r)$ (centered at the origin of $\cn$ with radius $r$) to $M$ satisfying the
following:
\begin{enumerate}
\item [(i)]$\Phi_t$ is a local biholomorphism from $D(r)$ to $M$;
\item[(ii)] $\Phi_t(0)=p$;
\item[(iii)] $\Phi_t^*(g(t))(0)=g_e$;
\item[(iv)] $\frac{1}{C}g_{e}\leq\Phi_t^*(g(t)) \leq  Cg_{e}$ in
$D(r)$. \item[(v)] for any $0<\alpha<1$, and $k\ge0$, the standard
$C^{k+\alpha}$ norm of $\Phi^*_t(g(t))$ in $D(r)$ is bounded by a
constant $C'$ which is independent of $t\ge -1$.
\end{enumerate}
where $g_e$ is the standard metric on $\cn$.
\end{prop}

\begin{rem}\label{s2r1}
Proposition 1.2 in \cite{TY} only requires the first covariant derivative of the
scalar curvature of $g(t)$ to be bounded independent of $t$.  In our case
however, we have bounds on all covariant derivatives of the Riemannain curvature
tensor independent of $k$.  Condition (iv) is derived by continuing the argument
in \cite{TY}, or \cite{CT3}. \end{rem}

As in \cite{CT3}, the following proposition is crucial:

\begin{prop}\label{s2p2}
Let $\lambda_1(t)\ge \dots>\lambda_n(t)>0$ be the eigenvalues of
$R_\ijb(p,t)$ relative to $g_\ijb(p,t)$.
\begin{enumerate}
\item[(i)] For any $\tau>0$,
$$\phi =\frac{\det(R_\ijb(p,t)+\tau g_{ij})}{\det(g_\ijb(p,t))}$$ is nondecreasing in
$t$. \item[(ii)] There is a constant $C>0$ such that
$\lambda_n(t)\ge C$ for all $t\ge0$.
\item[(iii)] For $1\le i\le
n$ the limit $ \lim_{t\to\infty}\lambda_i(t)$ exists.
\item[(iv)]
Let $\mu_1>\dots>\mu_l> 0$ be the distinct limits in (iii) and let
$\rho>0$ be such that the intervals $[\mu_k-\rho,\mu_k+\rho]$ for
$1\le k\le l$ are disjoint. For any $t$, let $E_k(t)$ be the sum
of the eigenspaces corresponding to the eigenvalues $\lambda_i(t)$
such that $\lambda_i(t)\in (\mu_k-\rho,\mu_k+\rho)$. Let $P_k(t)$
be the orthogonal projection (with respect to $g(t)$) onto
$E_k(t)$. Then there exists $T>0$ such that if $t>T$ and if $w\in
T_p^{(1,0)}(M)$, $|P_k(t)(w)|_t$ is continuous in $t$, where
$|\cdot|_t$ is the length measured with respect to the metric
$g(p,t)$.
\end{enumerate}
\end{prop}
\begin{proof}
The proof is identical to  the proof  of Proposition 3.1 in
\cite{CT3} for $\kappa=1$. Suppose $\kappa=0$. By Theorem 2.3 in
\cite{Cao}, if

\begin{equation}\label{s2e1}
Z_\ijb=\frac{\p R_\ijb}{\p t}+g^{k\bar l}R_{i\bar l}R_{k\bar j}
\end{equation}
then
\begin{equation}\label{s2e2}
Z_\ijb w^iw^{\bar j}\ge 0
\end{equation}
for any $w\in T^{(1,0)}(M)$. Let $p_{\ijb}=R_{\ijb}+ \tau
g_{\ijb}$ and  denote its inverse   by $(p^{\ijb})$. We have
\begin{equation}
\begin{split}
\frac{\p}{\p t}\log \phi&= p^{\ijb}\frac{\p}{\p t}p_{\ijb}-g^{\ijb}\frac{\p}{\p t}g_{\ijb}\\
&=p^{\ijb}\lf(\frac{\p}{\p t}R_{\ijb}-\tau R_\ijb \ri)+g^{\ijb} R_{\ijb} \\
&\ge  p^{\ijb}\lf(-g^{k\bar l}R_{i\bar l}R_{k\bar j}-\tau R_\ijb \ri)+g^{\ijb} R_{\ijb}\\
&=p^{\ijb}\lf(-g^{k\bar l}R_{i\bar l}R_{k\bar j}- \tau
p_\ijb\ri)+\tau^2 p^\ijb g_\ijb+g^\ijb R_\ijb
\end{split}
\end{equation}
  Now at the point $(p,t)$, we
choose a unitary basis such that $g_{\ijb}=\delta_{ij}$ and
$R_{\ijb}=\lambda_i\delta_{ij}$. Then $p_{\ijb}=(\lambda_i+ \tau
)\delta_{ij}$ and $p^{\ijb}=(\lambda_i+\tau)^{-1}\delta_{ij}$.
Hence we have
\begin{equation}
\begin{split}
\frac{\p}{\p t}\log \phi&\ge -\sum_{i=1}^n \frac{\lambda_i^2}{\lambda_i+\tau}- \tau n+\sum_{i=1}^n\frac{\tau^2}{\lambda_i+\tau}+\sum_{i=1}^n\lambda_i \\
&=\sum_{i=1}^n \lf(\frac{-\lambda_i^2}{\lambda_i+\tau}-\tau+\frac{\tau^2}{\lambda_i+\tau}+\lambda_i\ri)\\
&=0.
\end{split}
\end{equation}
From this (i) follows.  The proof of (ii)-(iv) is similar to the
proof of (ii)-(iv) in Proposition 3.1 of [CT3].
\end{proof}

For each $k$, consider the lifted family of metrics
$g_k(t):=\Phi_k ^* g(t_k+t)$ on $D(r)$ for $t\in [-1, \infty)$,
say. Then it is easy to see that $g_k(t)$ solves the \KRF
(\ref{s1e3}) on $D(r)$.  Then by Proposition \ref{s2p1} and the
\KRF it follows that some subsequence of $g_k(t)$ converges to a
smooth limit family $h(t)$, uniformly on compact subsets of
$D(r)\times(-1, \infty)$. It is easy to see that these are \K
metrics on $D$ for all $t$ and that $h(t)$ solves $(\ref{s1e3})$.
Moreover, by Proposition \ref{s2p2}, the eigenvalues of the Ricci
tensor $R_\ijb^h(t)$ of $h(t)$ at the origin are equal to
$\lim_{s\to\infty}\lambda_i(s)$ for any $t\in [0,\infty)$.
Therefore, $\mu_1>\mu_2>\cdots>\mu_l>0$ are distinct eigenvalues
of $R^{h}_\ijb(t)$ at the origin.  By the uniform bounds on the
covariant derivatives of the curvature tensor of $h(t)$ in
$D(r)\times(-1,\infty)$, and by Proposition \ref{s2p2}, we may
have the following inequality on $D(r)$, for $t\ge -1/2$, and by
choosing a smaller $r$ if necessary:
\begin{equation}\label{lowerbound}
    R_\ijb^h \ge Ch_\ijb
\end{equation}

\begin{thm}\label{s2t1}
Let $Rc^{h} _{i\jbar}(t)$ be the Ricci tensor of the metric $h_{
i\jbar}(t)$ on $D(r)$.
\begin{enumerate}
\item[(i)] For each $t\in [0, \infty)$ we have
$$Rc^{h} _{i\jbar}(t)+\kappa h_{ i\jbar}(t)=f_{i \jbar}(t)$$ for some smooth function
$f(t)$ on $D(r)$ such that $f_{ij}(t)=0$ and the gradient of
$f(t)$ in $h(t)$ is zero at the origin.
 \item[(ii)] Let
$\mu_1>\mu_2>\cdots>\mu_l>0$ be as above. For $1\le i\le l$, let
$E_i$ be the eigenspace corresponding to $\mu_i$ of $Ric^h(0,0)$
at $t=0$ at the origin with respect to $h(0)$. Then $E_i$ is also
the  eigenspace corresponding to $\mu_i$ of $Ric^h(0,t)$ for all
$t\ge0$ at the origin with respect to $h(t)$, $1\le i\le l$.
\end{enumerate}
\end{thm}

To prove the theorem, we first prove a lemma  which is a direct modification of
the results in \cite{Cao}. In the case of $\kappa=1$, it will be more convenient
to consider the transformed metric $\tilde{h}(t)=t h(\log t)$ which solves
(\ref{s1e3}) on $D(r)\times [e^{-1},\infty)$ with $\kappa=0$. It is clear that
$\tilde{h}(t)$ is the  limit of the transformed sequence $\tilde{g}_k (t):=t
g_k(\log t)$ uniform on compact sets of $D(r)\times [e^{-1},\infty)$ which also
satisfy (\ref{s1e3}) with $\kappa=0$.

Let $Z_{i\jbar}$ and $Z^k _{i\jbar}$ be the Harnack quadratic
tensors corresponding to $\tilde{h}(t)$ and $\tilde{g}_k(t)$
respectively as defined in Theorem 2.1 in \cite{Cao}. Namely for
any holomorphic vector $(V^i)$ at a point $q\in D(r)$,
\begin{equation}\label{Harnacktensors1}
  Z_\ijb=\frac{\p}{\p t}R^{\tilde h}_\ijb
  +{\tilde h}^{l\bar k}R^{\tilde h}_{i\bar k}R^{\tilde h}_{l\bar
  j}+R^{\tilde h}_{\ijb,k}V_{\bar k}+R^{\tilde h}_{\ijb,\bar k}V_k+
  R^{\tilde{h}}_{i\bar jk\bar
  l}V_{\bar k}V_l+\frac1t R^{\tilde h}_\ijb
\end{equation}
and $Z^k_\ijb$ is defined similarly. Denote the trace $\tilde
h^{\ijb}Z_\ijb$ of  $Z_\ijb$ by $Z$. Note that $Z$ is
a smooth function defined on the holomorphic tangent bundle
$T^{(1,0)}(D(r))$.

In case $\kappa=0$, then let $Q_{i\jbar}$ and $Q^k _{i\jbar}$ be
the Harnack quadratic tensors corresponding to $ {h}(t)$ and
${g}_k(t)$ respectively as defined in Theorem 2.3 in \cite{Cao}.
Namely for any holomorphic vector $(V^i)$ at a point $x\in D(r)$,
\begin{equation}\label{Harnacktensors2}
  Q_\ijb=\frac{\p}{\p t}R^{ h}_\ijb+h^{l\bar k}R^{  h}_{i\bar k}
  R^{  h}_{l\bar j}+R^h_{\ijb,k}V_{\bar k}+R^{ h}_{\ijb,\bar k}V_k
  +R^h_{i\bar
jk\bar
  l}V_{\bar k}V_l
\end{equation}
and $Q^k_\ijb$ is defined similarly. Denote the trace $
h^{\ijb}Q_\ijb$ of $Q_\ijb$ by $Q$.

\begin{lem}\label{s2l2}\ \

 \begin{itemize}

    \item [(i)] For any holomorphic vector $V\in T^{(1,0)}(D(r))$,
    $Z_\ijb$ is a nonnegative quadratic form. Moreover, if $Z$
    is positive  at some point $q$ for all $V\in
    T_{q}^{(1,0)}(D(r))$ at $t=t_0$, then  $Z$ is positive for
    all $t>t_0$ and for  $V\in T^{(1,0)}(D(r))$.

    \item [(ii)] For any holomorphic vector $V\in T^{(1,0)}(D(r))$,
    $Q^{ h}_\ijb$ is a nonnegative quadratic form. Moreover, if $Q^{  h}$
    is positive  at some point $x_0$ for all $V\in
    T_{q}^{(1,0)}(D(r))$ at $t=t_0$, then  $Q^{  h}$ is positive for
    all $t>t_0$ and for  $V\in T^{(1,0)}(D(r))$.
 \end{itemize}
\end{lem}

\begin{proof} For any holomorphic vector $W$, $Z^k_\ijb W^iW^{\bar
j} \geq 0$ for all $k$ by Theorem 2.1 in \cite{Cao}. Since
$Z_\ijb$ is the limit of the $Z^k_\ijb$'s on $D(r)$ for all $t$,
  $Z_\ijb W^iW^{\bar j} \geq 0$. This proves the first
statement of (i). The first statement of (ii) can be proved
similarly.

To prove the second statement in (i),   assume there is some
$x_0\in D(r)$ and $t_0 \geq e^{-1}$ so that $Z>0$ for all $V\in
T_{x_0}^{(1,0)}(D(r))$.  Given any $T>t_0$, we note that for $C>0$
there exists some $K>0$ such that given any point $(x, t) \in
D(r)\times [t_0,T]$ and $V\in T_x ^{(1, 0)} (D(r))$, with
Euclidean length $||V||>K$, we must have

\begin{equation}\label{zbound}
 Z>C
\end{equation} at $(x, t)$ and $V$.  This follows from (\ref{Harnacktensors1}), (\ref{lowerbound}) and the fact that the
curvature tensor of $\tilde{h}(t)$ and its covariant derivatives
in time and space are uniformly bounded on $D\times [t_0,T]$ by
constants independent of space and time by our estimates on the
$\tilde{g}_k 's$.  Hence there exist a neighborhood $U$ of $x_0$
and $\e>0$ such that $Z \ge \e$ for all holomorphic vector $V$ at
$x\in U$ at $t=t_0$.

 Choose a smooth function $F$ on $D(r)$ such that
$F(x_0)>0$, $F$ is zero outside a small neighborhood of $x_0$, and
$Z-\frac{F}{t_0 ^{2}}\geq 0$ for all $V$ in  $T^{(1,0)}(D(r))$ at
$t=t_0$.

 Let $F$ evolve by the   heat equation on $D\times
[t_0, T]$ with the following initial and boundary conditions:
\begin{equation}
\begin{split}
\frac{\p}{\p t}F&=\Delta _t F \hspace{11pt} in \hspace{11pt}
 D\times [t_0, T]\\
F &=0 \hspace{11pt} on \hspace{11pt} \partial D \times [t_0, T]\\
F(x,t_0)&=F(x)
\end{split}
\end{equation}
where $\Delta _t$ is the Laplacian relative to $\tilde{h}(t)$. $F$
is then strictly positive in $D\times (t_0, T]$ by the strong
maximum principle \cite[Theorem 5, Chapter 3]{PW}. We will show
that $\tilde{Z}:=Z-\frac{F}{t^2}$ is also nonnegative for all $V$
and $(x, t) \in D\times (t_0, T]$. Without loss of generality we
may assume $\tilde h$ is smooth up to the boundary of $D(r)$. Let
$\tilde{Z}$ assume its minimum over all $(x, t)\in \ol D(r)\times
[t_0, T]$ and $V$, at some point $(x, t_1)$ and some vector $V_0$.
This minimum exists by compactness and by (\ref{zbound}). Now
assume this minimum is negative.  Then by the initial condition of
$F$, we must have $t_1>t_0$. Also, by the non-negativity of $Z$
and the zero Dirichlet boundary condition on $F$, we must have $x$
strictly inside $D(r)$.  As in \cite{Cao}, we may then extend
$V_0$ locally around $(x, t_1)$ in space-time such that (3.2) in
\cite{Cao}  holds. At $(x, t_1)$ we then have
\begin{equation}
\begin{split}
(\frac{\p }{\p t}-\Delta _t)\tilde{Z}&=(\frac{\p  }{\p t}-\Delta _t)Z+2\frac{F}{t^3}\\
&=Z_{i\jbar}R^{\tilde{h}} _{j\ibar}-2\frac{Z}{t}+2\frac{F}{t^3}\\
&=Z_{i\jbar}R^{\tilde{h}} _{j\ibar}-2\frac{\tilde{Z}}{t}\\
&>0.
\end{split}
\end{equation}
But this contradicts the fact that $\tilde{Z}$ is minimal at $(x,
t_1)$ and $V_0$.  Thus $\tilde{Z}$ is non-negative as claimed
which completes the proof of (i) of the lemma.

We now consider the second statement of (ii). Assume there is some
$x_0\in D(r)$ and $t_0 \geq 0$ such that $Q$ is positive at $(x_0,
t_0)$ for all $V$.  As before, for any $T>t_0$ we observe that
given any $C>0$, there exists some $K>0$ such that given any point
$(x, t) \in D(r)\times [t_0,T]$ and $V\in T_x ^{(1, 0)} (D(r))$
with Euclidean length $||V||>K$, we must have $Z>C$ at $(x, t)$
and $V$.

Hence we can choose a smooth function $F$ on $D$ such that
$F(x_0)>0$, $F$ is zero outside a small neighborhood of $x_0$, and
$Q-F\geq 0$ for all $V$ everywhere on $D(r)$ at $t=t_0$. Let $F$
evolve by heat equation on $D\times [t_0, T]$ with the following
initial and boundary conditions:
\begin{equation}
\begin{split}
\frac{\p}{\p t}F&=\Delta _t F \hspace{11pt} in \hspace{11pt} D(r)\times [t_0, t]\\
F &=0 \hspace{11pt} on \hspace{11pt} \partial D(r) \times [t_0,
T]\\
F(x,t_0)&=F(x)
\end{split}
\end{equation}
where $\Delta _t$ is the Laplacian relative to $h(t)$. $F$ is then
strictly positive in $D\times (t_0, T]$ by the strong maximum
principle.  Now given any $\epsilon >0$, we will show that
$\tilde{Q}:=Q-F+\epsilon e^{t}$ is nonnegative for all $V$ and
$(x, t) \in D(r)\times (t_0, T]$. Letting $\epsilon$ approach
zero, this will prove that $Q-F$ is nonnegative for all $V$ and
$(x, t) \in D(r)\times (t_0, T]$, thus proving the lemma.

Let $\tilde{Q}$ assume its minimum over all $(x, t)\in D(r)\times
[t_0, T]$ and $V$, at some point $(x, t)$ and some vector $V_0$.
Now assume this minimum is negative. Then by our initial condition
of $F$, we must have $t>t_0$. Also, by the non-negativity of $Q$
and the zero Dirichlet boundary condition on $F$, we must have $x$
strictly inside $D(r)$. We may then extend $V_0$ locally around
$(x, t)$ in space-time such that (3.1) in \cite{Cao} holds.  At
$(x, t)$ we then have
\begin{equation}
\begin{split}
(\frac{\p }{\p t}-\Delta _t)\tilde{Q}&=(\frac{\p }{\p t}-\Delta _t)Q +\epsilon e^t\\
&=Q_{i\jbar}R^{\tilde{h}}_{j\bar i} +\epsilon e^t\\
&>0.
\end{split}
\end{equation}
But this contradicts the fact that $\tilde{Q}$ is minimal at $(x,
t)$ and $V_0$. Thus $\tilde{Q}$ is non-negative as claimed which
completes the proof of (ii).
\end{proof}
\begin{proof}[Proof of Theorem \ref{s2t1}]
(i): We begin with the case where $\kappa=1$ in Theorem
\ref{s2t1}.  With the same notations as in the Lemma \ref{s2l2},
by Proposition \ref{s2p2} and the arguments following (3.7) in
\cite{CT3},  we see that
$$t R^{\tilde{h}} (0, t)$$ is constant for
all $t\in [e^{-1}, \infty) $ where $R^{\tilde{h}}$ is the scalar
curvature of $\tilde{h}$. Thus at the space time point $(0, t)$ we
have $R^{\tilde{h}}+t\frac{\p}{\p t}R^{\tilde{h}}=0$. Applying
Lemma \ref{s2l2} and following the exact argument in the proof of
Theorem 4.2 in \cite{Cao}, we conclude that for each $t\in
[e^{-1}, \infty)$ there is a smooth function $\tilde{f}(t)$ on
$D(r)$ such that the gradient of $\tilde{f}(t)$ is holomorphic and
is zero at the origin. Moreover, we have
\begin{equation}\label{s2e9}
Rc^{\tilde{h}}_{i\jbar}=\tilde{f}_{i\jbar}+\frac{1}{t}\tilde{h}_{i\jbar}.
\end{equation}
Transforming $\tilde{h}$ back to $h$, it is easy to see that (i)
in Theorem \ref{s2t1} is true for $\kappa=1$.

We now consider the case of $\kappa=0$.    By Proposition
\ref{s2p2} we have $R^{h}(0, t)$ is constant for $t\in [0,
\infty)$, and in particular, $\frac{\p}{\p t}R^{h}=0$ at the space
time point $(0, t)$.  Applying Lemma \ref{s2l2} and following the
exact argument in the proof of Theorem 4.1 in \cite{Cao}, we
conclude that for each $t\in [0, \infty)$ there is a smooth
function $ {f}(t)$ on $D(r)$ such that the gradient of $ {f}(t)$
is holomorphic and is zero at the origin. Moreover,   we have
\begin{equation}\label{s2e99}
Rc^{h}_{i\jbar}= {f}_{i\jbar}.
\end{equation}
This completes the proof of (i) in Theorem \ref{s2t1}.

 (ii): Let $\lambda_1\ge\dots\ge \lambda_n>0$ be the eigenvalues
 of $Rc^h$ at $0$. Note that they are independent of $t$.  Suppose
  $v\in T_0^{1,0} (D(r))$ is not an
eigenvector for $Rc^h(0, 0)$ for $\lambda_k$. We will show that
$v$ can not be an eigenvector for $Rc^{h}(0, t)$ for all $t\in (0,
\infty)$ for $\lambda_k$. It is sufficient to prove that, the
quantity
\begin{eqnarray}
\begin{split}
  F(t):&= |Rc^{h} (0, t)(v, \cdot)-\lambda_k h(0, t)(v, \cdot)|^2_{h(t)}  \\
    &=  h^{j\bar k}(0, t) \ol{(R^h_\ijb (0, t) -\lambda_k h_\ijb(0, t)
    )
    v^i}\\
     &\ \ \ \cdot{(R^h_{l\bar k}(0, t)   -\lambda_k h_{l\bar k}(0, t) )
     v^l}
     \end{split}
\end{eqnarray}
  can never be zero, since  this is zero at $t$ if and only if
$v$ is an eigenvector for $Rc^{h}(0, t)$ with eigenvalue
$\lambda_k$.  Let $t_0\in [0, \infty)$ be arbitrary. Choose a
holomorphic coordinate in $D(r)$ such that at $0\in D(r)$ we have
$h_{ i\jbar}(t_0)=\delta _{i\jbar}$ and $R^{h}
_{i\jbar}(t_0)=\lambda_i \delta_{i\jbar}$.
 Let $f(t)$ be as in Lemma \ref{s2l2}. It is not hard to show
that we may choose some $1>\delta>0$ such that starting at any
point $p \in  D(\delta r)$ we may flow along $-\frac 12\nabla
f(t_0)$ for $t\in [0, 1]$ while staying inside $D(r)$ and let
$\varphi_t$ be the local biholomorphism determined by the flow.
Note that the origin is a fixed point of the flow because $\nabla
f(t)=0$ at the origin.  Let $g(t)= \varphi_t^*(h(t_0))$ be the
soliton metric on $ D(\delta r)\times [0, 1)$ with initial
condition $g(0)=h(t_0)$. Then at $0 \in D(\delta r)$, in the above
coordinates, we have $g_{i\jbar}(t)=e^{-(\lambda_i+\kappa)
t}\delta_{i\jbar}$ and $R^g_{i\jbar}(t)=\lambda_i
e^{-(\lambda_i+\kappa) t}\delta_{i\jbar}$. For any $k$ we then
have
\begin{equation}
\begin{split}
G(t)&:=|Rc^g(0, t)(v,\cdot)  - \lambda_k g(0, t)(v,\cdot)|_{g(t)} ^2\\
&=\sum_{i\neq k} (\lambda_i -\lambda_k)^2 |v^i|^2e^{-(\lambda_i+\kappa) t}\\
\end{split}
\end{equation}
and thus
\begin{equation}
\begin{split}
G'(0)&=\sum_{i\neq k} -(\lambda_i+\kappa)
(\lambda_i -\lambda_k)^2 |v^i|^2\\
&\geq- (\lambda_1+\kappa) G(0).
\end{split}
\end{equation}
From the above equation, the fact that $g(0)=h(t_0)$ on $ D(\delta r)$, and the
fact that both $g(t)$ and $h(t_0+t))$ both solve (\ref{s1e3}) on $ D(\delta r)
\times [0, 1)$, it follows that $G(0)=F(t_0)$ and $G'(0)=F'(t_0)$. Hence for any
choice of $k$ we have $F'(t_0)\geq -(\lambda_1+\kappa) F(t_0)$. But $t_0\in [0,
\infty)$ is arbitrary. Thus for any $k$ we have $F'(t)\geq - (\lambda_1+\kappa)
F(t)$ for all $t\in [0, \infty)$.  It is now easy to see that if $F(0)\neq 0$,
then $F(t)$ cannot be zero for any $t$. This completes the proof of our claim.

From the claim, we conclude that if $v$ is an eigenvector of
$h(t)$ for $t>0$ with eigenvalue $\lambda_k$, then $v$ must be in
$E_k$. Since the multiplicity of each eigenvalue $\mu_k$ is
constant in $t$, from this it is easy to see the theorem is true.
\end{proof}

 Now
given $(M, g(t))$   as in Theorem \ref{s2t1}, we denote the
eigenvalues of $Rc(p,t)$ by $\lambda_i (t)$ for $i=1,...,n$ as
before, and we let $\mu_k$, $E_k(t)$ and $P_k(t)$ for $k=1,...,l$
be as in Proposition \ref{s2p2}. We let $n_m$ for $m=0,...,l-1$ be
such that $\lambda_k (t) \in (\mu_m -\rho, \mu_m+\rho)$ for all
$n_{m}\leq k\leq n_{m+1}$ and $t$ sufficiently large  such that
the intervals $[\mu_m
  -\rho, \mu_m+\rho]$ are disjoint as in Proposition \ref{s2p2} part (iv).
  For any nonzero vector $v \in T_p ^{1, 0}(M)$, let $v(t)=v/ |v|_t $
  where $|v|_t$ is the length of $v$ with respect to $g(t)$
  and  $v_i(t)=P_i(t)v(t)$.

We now show that Theorem 4.1 in \cite{CT3} is also true for $(M,
g(t))$ in our case; that $Rc(p,t)$ can be `diagonalized'
simultaneously near infinity and that $g(t)$ is `Lyapunov
regular', to borrow a notion from dynamical systems (see
\cite{BP}).

\begin{thm}\label{s2t3}
Let $(M,g(t))$ be as described in the beginning of the section.  Then
   $V=T_p^{(1,0)}(M)$ can be decomposed orthogonally with respect to $g(0)$ as
 $V_1\oplus \cdot\cdot\cdot\oplus V_l$ so that the following are true:
\begin{enumerate}
\item[(i)] If $v$  is a nonzero vector in  $V_i$ for some $1\le i\le l$, then
 $\lim_{t \to \infty} |v_i(t)|=1$ and thus $\lim_{t\to\infty}Rc(v(t),\bar v(t))=\mu_i$ and
$$\lim_{t \to \infty} \frac{1}{t}\log \frac{|v|_t^2}{|v|_0^2}=-\mu_i-\kappa.$$  Moreover, the convergence are uniform over all
 $v\in V_i\setminus\{0\}$.
\item[(ii)] For $1\le i, j\le l$ and for nonzero vectors  $v \in V_i$ and $w \in
V_j$ where $i\neq j$, $\lim_{t\to \infty}\langle v(t),w(t)\rangle_t=0$ and the
convergence is uniform over all such nonzero vectors $v, w$. \item[(iii)]
$\dim_\C(V_i)=n_{i}-n_{i-1}$  for each $i$. \item[(iv)]
$$\sum_{i=1}^l(-\mu_i-\kappa)\dim_\C V_i=\lim_{t\to\infty} \frac1t\log
\frac{\det(g_{i\bar j}(t))}{\det ({g}_{i\bar j}(0)}.$$
\end{enumerate}
\end{thm}
\begin{proof} Let $t_k\to\infty$ and construct $g_k$ with limit
$h(t)$ as in Theorem \ref{s2t1}. Observe that since $h(t)$ is a
smooth limit of the $g_k(t)'s$ on $D(r)\times[0, \infty)$, the
analog of Lemma 3.2 in \cite{CT3} is true in our case.  Using this
and (ii) in Theorem  \ref{s2t1}, we may prove the theorem by
contradiction exactly as in the proof of Theorem 4.1 in
\cite{CT3}.
\end{proof}

\section{Transition maps}

 Let $(M,g(t))$ be as in Theorem \ref{s1t3}  and let $p\in M$ be fixed.
 In addition, we will assume that $R_\ijb\ge a'g_\ijb$ for some $a'>0$
 at $p$ and $t=0$. Then by Proposition \ref{s2p2}, there exists $a >0$ such that
\begin{equation}\label{riccibound}
  R_\ijb\ge ag_\ijb
\end{equation}
 at $p$ for all $t\ge0$. Since the covariant derivatives Riemannian
  curvature tensor  $Rm_t$ of $g(t)$ is uniformly bounded, we can conclude
  that by choosing a possibly smaller $a>0$,
(\ref{riccibound}) is still true  in $B_t(p,R)$ for some $R>0$
independent of $t$, where $B_t(p,R)$ is the geodesic ball of
radius $R$ with center at $p$ with respect to the metric $g(t)$.
Since $B_{t_1}(p,R)\subset B_{t_2}(p,R)$ for $t_2\ge t_1$ as
$g(t)$ is shrinking, we have
\begin{equation}\label{lengthdecay}
   L_{t_2}(\gamma)\le e^{-\frac a2(t_2-t_1)}L_{t_1}(\gamma)
\end{equation}
for any curve in $B_{t_1}(p,R)$. Here $L_t$ denotes the length
function with respect to $g(t)$.

  Recall that by Proposition \ref{s2p1},   there exist $r>0$ and $C>0$
  independent of $t\ge0$ and
  a holomorphic maps $\Phi_t:D(r)\to M$ with the following
  properties:
  \begin{enumerate}
\item [(i)] $\Phi_t$ is a local biholomorphism from $D(r)$ onto its image;
\item[(ii)] $\Phi_t(0)=p$; \item[(iii)] $\Phi_t^*(g(t))(0)=g_e$; \item[(iv)]
$\frac{1}{C}g_{e}\leq\Phi_t^*(g(t)) \leq  Cg_{e}$ in $D(r)$.
\end{enumerate}
where $g_e$ is the standard metric on $\cn$. Let $T>0$ and denote
$\Phi_{iT}$ simply by $\Phi_i$. In this section, we want to
construct  injective holomorphic maps $F_i$ from $D(\rho)$ to
$D(\rho)$ for some $\rho$ such that $\Phi_i=\Phi_{i+1}\circ
F_{i+1}$. We should emphasis that $\Phi_i$ may not be a covering
map.

In this section, we always assume that $t\ge0$. By property (iv)
and reducing $r$ if necessary, we may assume that
$\Phi_t(D(r))\subset B_t(p,R)$, where $R>0$ is such that
(\ref{lengthdecay}) is true. In fact, we have the following:

\begin{lem}\label{s3l1}
For any $0<\rho<r$, there exists $R_1>0 $ independent of $t$  such that
\begin{equation}\label{image1}
 B_{t}(p,\frac1R_1)\subset\Phi_t(D(\rho))\subset B_{t}(p, R_1)
\end{equation}
\end{lem}
\begin{proof}    By (iv) above,  it is easy to see that
$$
\Phi_t(D(\rho))\subset B_{t}(p, C_1)
$$
for some $C_1>0$ independent of $t$, where $B_t(p,r)$ is
 the geodesic ball with radius $r$ centered at $p$ with respect to
 $g(t)$. On the other hand, $\wh B_{t}(0,C_2)\subset D(\rho)$
for some $C_2>0$ independent of $t$, where $\wh B_{t}(0,C_2)$ is
 the geodesic ball with radius $C_2$ centered at $0$ with respect to
 $\Phi_t^*(g(t))$. Hence
 $$\Phi_t(D(\rho))\supset\Phi_t(\wh B_{t}(0,C_2))\supset
 B_{t}(p,C_2).
 $$
 From this it is easy to see the lemma follows.
\end{proof}

\begin{lem}\label{s3l2}
For any $0<\rho\le r$, where $r$ is as in (i)-(iv),  there exists $\rho_1>0$
independent of $t$, satisfying the following for any $t$: Let $\gamma$ be a
smooth curve in $M$ with $\gamma(0)=q$, such that $q\in B_{t}(p,\rho_1)$ and
$L_t(\gamma)<\rho_1$. Then $\Phi_t(z)=q$ for some $z\in D(\frac\rho8)$, and for
all such $z$ there is a unique lift $\wt\gamma$ of $\gamma$ by $\Phi_t$ so that
$\wt\gamma(0)=z$ and $\wt\gamma\subset D(\frac\rho2)$.
\end{lem}
\begin{proof} Let $\rho_1>0$ be determined later.   Let $q\in B_{t}(p,\rho_1)$
and let $\gamma(s)$, $0\le s\le \rho_1$ be a curve from $q$
parametrized by arc-length with respect to $g (t)$. Suppose $z\in
D(\frac\rho8)$ with $\Phi_t(z)=q$

 Since
$\Phi_t$ is a local biholomorphism, there exists $s_0>0$ and a
curve $\wt \gamma$ from $z$ with  $\wt\gamma \subset
D(\frac18\rho)$ such that $\Phi_t\circ\wt\gamma=\gamma$ on
$[0,s_0]$. Let $A$ be the set of $s$, such that $\gamma$ has a
lift $\wt\gamma$ in $D(\frac12\rho)$ on $[0,s]$ with
$\wt\gamma(0)=z$ . Since $\Phi_t$ is a local biholomorphism, $A$
is open in $[0,\rho_1]$. Suppose $s_k\to s$, and $s_k\in A$. By
(iv), there is a constant $C>0$ which is independent of $t$   such
that
$$
C^{-1}L_e(\wt\gamma|_{[0,s_k]})\le L_t(\gamma|_{[0,s_k]})\le
\rho_1.
$$
where $L_e$ is the length with respect to $g_e$. Hence
$$
L_e(\wt\gamma|_{[0,s_k]})\le C\rho_1\le \frac14\rho
$$
if $\rho_1<\frac1{4 C}\rho$. Note that since $\wt\gamma(0)=z\in D(\frac18\rho)$,
we may assume that $\wt\gamma(s_k)\to x$ for some $z_1\in \ol {D(\frac38\rho)}$.
From this it is easy to see that $\gamma$ can be lifted up to $s$ while staying
in $D(\frac12\rho)$ . Hence $A$ is also closed. Since $\Phi_t$ is a local
biholomorphism, the lifting must be unique. In particular, by choose a smaller
$\rho_1$ which is independent of $t$, we conclude that every minimal geodesic
from $p$ with length less than $\rho_1$ can be lifted to a curve in
$D(\frac18\rho)$. Hence for all $q\in B_t(p,\rho_1)$, there is a point $z\in
D(\frac18\rho)$ such that $\Phi_t(z)=q$. This completes the proof of the lemma.
\end{proof}

\begin{lem}\label{s3l3} Fix $t\ge0$.  Let $0<\rho\le r$ be given and let $\rho_1$ be as in Lemma \ref{s3l2}.  Given any $\e>0$,
 there exists $\delta>0$ satisfying the following properties:

Let  $\gamma(\tau)$, $\beta(\tau)$, $\tau\in [0,1]$ be  smooth
curves from $q\in B_{t}(p,\rho_1)$ with length less than $\rho_1$
and let $z_0\in D(\frac18\rho)$ with $\Phi_t(z_0)=q$. Let
$\wt\gamma$, $\wt\beta$ be the liftings from $z_0$ of $\gamma$ and
$\beta$ as described in Lemma \ref{s3l2}.
 Suppose  $d_t(\gamma(\tau),\beta(\tau))<\delta$ for all $\tau\in [0,1]$,
 then $d_e(\wt\gamma(1),\wt\beta(1))<\e$.
 Here   $d_t$ is the distance in $g(t)$ and $d_e$ is the Euclidean distance.

\end{lem}
\begin{proof} Since $\Phi_t$ is a local biholomorphism, there is
$\sigma>0$ such that $\Phi_t$ is a biholomorphism onto its image
when restricted on $D (z,\sigma)$ for all $z\in
\ol{D(\frac12\rho)}$, where $D(z,\sigma)$ is the Euclidean ball
with center at $z$ and radius $\sigma$.

Let $q$, $z_0$, $\gamma$, $\beta$, $\wt\gamma$ and $\wt\beta$ as
in the lemma. Since $\wt\gamma\subset D(\frac12\rho)$, by property
(iv) of $\Phi_t$, there exists $C_1>0$  such that
\begin{equation}\label{localbiholo1}
\Phi_t(D(\wt\gamma(\tau)),\sigma)\supset
B_{t}(\gamma(\tau),C_1^{-1}\sigma)
\end{equation}
and
 \begin{equation}\label{localbiholo2}
d_e (\wt\gamma(\tau),z)\le C_1d_t(\gamma(\tau),\Phi_t(z))
\end{equation}
for all $z\in D(\wt\gamma(\tau),\sigma)$ with $\Phi_t(z)\in
B_{t}(\gamma(\tau),C_1^{-1}\sigma)$.   Note that  $C_1$ is
  independent of the curves $\gamma$ and $\beta$.

Given $0<\e<\sigma$,  let  $0<\delta<C_1^{-1} \e<C_1^{-1}\sigma$.
Suppose $\beta$ and $\wt\beta$ are as in the lemma such that
$d_t(\gamma(\tau),\beta(\tau))<\delta$ for all $\tau$. Since
$\wt\gamma(0)=\wt\beta(0)= z_0 $, we have
$d_e(\wt\gamma(\tau),\wt\beta(\tau))<\e$ in $[0,\tau_0]$ for some
$\tau_0>0$. Let $A$ be the set $\tau$ in $[0,1]$ such that
$d_e(\wt\gamma(\tau'),\wt\beta(\tau'))<\e$ for all
$\tau'\in[0,\tau]$. Then $A$ is nonempty and is open. Suppose
$\tau_k\in A$ and $\tau_k\to \tau$. Then
$$
d_e(\wt\gamma(\tau),\wt\beta(\tau))\le \e.
$$
Since $\e<\sigma$, $\delta<C_1^{-1}\sigma$,   by
  (\ref{localbiholo2})  and the fact that
$\Phi_t(\wt\gamma(\tau))=\gamma(\tau)$,
$\Phi_t(\wt\beta(\tau))=\beta(\tau)$, we have

$$
d_e(\wt\gamma(\tau),\wt\beta(\tau))\le
C_1d_t(\gamma(\tau),\beta(\tau))\le C_1\delta<\e.
$$
Hence $\tau\in A$ and $A=[0,1]$. This completes the proof of the lemma.
\end{proof}
Apply Lemma \ref{s3l2} to $\rho=r$ and choose $\rho_1$ as in the
lemma. Note that $\rho_1$ is   independent of $i$ and $T$. For any
$z\in D(r)$, let $\gamma^*(\tau)$, $0\le \tau\le 1$, be the line
segment from 0 to $z$, and let $\gamma=\Phi_i\circ\gamma^* $. By
(\ref{lengthdecay}) and property (iv) for $\Phi_t$, there is a
constant $C_1>0$ independent of $i$ and $T$ such that
\begin{equation}
  L_{i+1}(\gamma)\le e^{-\frac a2T}L_i(\gamma)\le C_1e^{-\frac
  a2T}r.
  \end{equation}
Now we choose $T>0$ large enough so that $C_1e^{-\frac
  a2T}r<\rho_1$. Then by Lemma \ref{s3l2}, there is a unique lift
  $\wt\gamma$
  of $\gamma$ by $\Phi_{i+1}$ so that $\wt\gamma(0)=0$ and
  $\wt\gamma\subset D(\frac r2)$. We then define
  $F_{i+1}(z)=\wt\gamma(1)$. $F_{i+1}$ is then well-defined by
    the uniqueness of  the lifting. We have:
  \begin{lem}\label{transitionmaps} The maps $F_{i+1}$ satisfy the
  following:

  \begin{enumerate}
    \item [(a)] $F_{i+1}:D(r)\to D(\frac r2)$, $F_{i+1}(0)=0$ and
    $\Phi_{i}=\Phi_{i+1}\circ F_{i+1}$.
    \item [(b)] For each $i$, $F_{i+1}$ is a local  biholomorphism.
    \item [(c)]
    $$
    b_1|v|\le |F'_{i+1}(0)v|\le b_2|v|
    $$
    for some $0<b_1\le b_2<1$ for all $i$ and for all vector $v\in
    \cn$, where $F'$ is the Jacobian of $F$.
    \item[(d)] There exist  $r>r_1>0$ and $0<\theta<1$ independent of $i$ such that
    $$
    |F_{i+1}(z)|\le \theta|z|
    $$
    for all $i$ and for all $z\in D(r_1)$ and $F_{i+1}$ is injective on $D(r_1)$.
  \end{enumerate}
  \end{lem}
  \begin{proof} (a) follows immediately from the definition of
  $F_{i+1}$.

  To prove (b), let us first prove that $F_{i+1}$ is continuous. Let $z_0\in D(r)$, $\gamma^*(\tau)$, $0\le \tau\le1$, be the line segment from 0 to $z_0$,
  $\gamma=\Phi_i\circ\gamma^*$ and $\wt\gamma$ is the lift of $\gamma$ by
  $\Phi_{i+1}$ with $\wt\gamma(0)=0$.   Let $w=\wt\gamma(1)=F_{i+1}(z_0)$.
  Given $\e>0$, let $\delta>0$ be  as in Lemma \ref{s3l3}   for
  $\Phi_{i+1}$.

 We may assume that $|z_0|\le 1-\eta$ for some $\eta>0$. Since
 $\Phi_i$ is uniformly continuous on $D(1-\eta/2)$,  there exists $\sigma>0$ such that if $|z_1-z_2|<\sigma$, $z_1,z_2\in D(1-\eta/2)$, then
$d_i(\Phi_i(z_1),\Phi_i(z_2))<\delta$.

Moreover, it is easy to see that we can find $\delta'>0$ such that
if $|z_0-\zeta|<\delta'$, then the ray $\beta^*$ defined on
$[0,1]$ from $0$ to $\zeta$ satisfies
$|\gamma^*(\tau)-\beta^*(\tau)|<\sigma$ and $\beta^*\subset
D(1-\eta/2)$. Let $\zeta$ be such a point in $D(r)$ with $\beta^*$
as above and let $\beta=\Phi_i\circ\beta^* $. Hence we have
$$
 d_{i+1}( \gamma(\tau),\beta(\tau))\le d_{i}( \gamma(\tau),\beta(\tau))<\delta.
$$
By Lemma \ref{s3l3}, if $\wt\beta$ is the lift of $\beta$ by
$\Phi_{i+1}$ with $\wt\beta(0)=0$, then
$|\wt\gamma(1)-\wt\beta(1)|<\e$. That is to say,
$|F_{i+1}(z_0)-F_{i+1}(\zeta)|<\e$ and $F_{i+1}$ is continuous.

Now it is easy to see that $F_{i+1}$ is a local biholomorphism. In
fact, suppose $F_{i+1}(z)=w$  and suppose $\Phi_i(z)=x$ and
$\Phi_{i+1}(w)=y$. Let $\e_1>0$ be such that $\Phi_i$ and
$\Phi_{i+1}$ are biholomorphisms when restricted on $D(z,\e_1)$
and $D(w,\e_1)$ respectively. Since $F_{i+1}$ is continuous, we
can find $0<\delta_1<\e_1$, such that
$F_{i+1}(D(z,\delta_1))\subset D(w,\e_1). $ Since
$\Phi_i=\Phi_{i+1}\circ F_{i+1}$, we have
$$
F_{i+1}=\Phi_{i+1}^{-1}\circ\Phi_i
$$
on $B_e(z,\delta_1)$. Hence $F_{i+1}$ is a local biholomorphism.

To prove (c), by properties   (ii), (iii) of $\Phi_t$, the fact that
$F_{i+1}=\Phi_{i+1}^{-1}\circ\Phi_i$ near the origin,   the fact that $R_\ijb(t)$
is uniformly bounded and (\ref{riccibound}), it is easy to see that (c) is true.

 To prove (d), by  gradient
estimates, we have
$$
F_{i+1}(z)=F_{i+1}(0)+F_{i+1}'(0)z+H_i(z)=F_{i+1}'(0)z+H_i(z)
$$
where $|H_i'(z)|\le C|z|$ for some constant $C$ independent of $i$ on
$D(\frac12r)$, say. Hence by (c), there exist  $r>r_1>0$ and $1>\eta>0$,
independent of $i$ such that $F_{i+1}:D(r_1)\to D(r_1)$ so that
$$ |F_{i+1}(z)|\le \theta|z|  $$
     for all $i$ and for all $z\in D(r_1)$ and $F_{i+1}$ is
     injective on $D(r_1)$.
  \end{proof}

  By this lemma and Theorem \ref{s2t3}, using the method in \cite[\S5]{CT3},
   we can prove the following:
  \begin{lem}\label{correctionmaps} Let $F_i$ be as in Lemma
  \ref{transitionmaps}. There exist biholomorphisms $G_i$ on $\cn$
  and polynomial maps $T_i$ with the following properties:

  \begin{enumerate}
    \item [(a)] $T_i(0)=0$, $T_i'(0)=Id$, and $\sup_{z\in
    D(1)}|T_i(z)|\le C_1$ for some constant $C_1$ independent of
    $i$.
    \item [(b)] $G_i(0)=0$ and for all open set $U$ containing the
    origin
    $$
    \bigcup_{k=1}^\infty (G_1\circ \cdots\circ G_k)^{-1}(U)=\cn.
    $$
    \item[(c)] There exist $0<r_2<r_1$ and
  $C_2>0$ independent of $k\ge1$ such that  $$
  G_{k+1}^{-1}\circ G_{k+2}^{-1}\circ \cdots
  \circ G_{k+l}^{-1}\circ T_{k+l}\circ F_{k+l}\circ\cdots \circ F_{k+2}\circ F_{k+1}
  $$
  converges uniformly on $D(r_2)$ as $l\to\infty$ to an injective holomorphic map
  $\Psi_k$ such that
  $$
  D(C_2^{-1}r_2)\subset \Psi_k(D(r_2))\subset D(C_2r_2).
  $$
  \item [(d)] By choosing $r_2$ smaller if necessary, we may have that $T_i$ is injective on $D(r_2)$ such that $T_i^{-1}$
  is defined on $D(r_2)$ and $T_i^{-1}(D(r_2)) \subset D(r)$ for all $i$.
  \end{enumerate}
  \end{lem}

\section{Proof of the main theorem }

We are ready to prove  Theorem  \ref{s1t3}. By Remark
\ref{reduction}, it is sufficient to assume that $M$ is simply
connected and $g(t)$ has positive Ricci curvature for all $t$.

Let $\Phi_i$ be as in the previous section so that one can define
$F_i$ which satisfy Lemma \ref{transitionmaps}. Let $G_i$, $T_i$,
  $r_2$, $C_1$ and $C_2$ be as in Lemma \ref{correctionmaps}.
   We want to construct a
biholomorphism from $\cn$ onto $M$ as follows: Let
$\Omega_i=(G_i\circ\cdots\circ G_1)^{-1}(D(r_2))$ and define
$$
S_i=\Phi_i\circ T_i^{-1}\circ G_i\circ \cdots \circ G_1
$$
which is defined on $\Omega_i$ by Lemma \ref{correctionmaps}.

  Theorem \ref{s1t3} will be proved if we can prove that
  $S_i$ converges to a bihololomorphism from $\cn$ onto $M$. We will
  prove this in several steps as described in the following
  lemmas.
\begin{lem}\label{limit} For all $z\in \cn$, $\lim_{i\to\infty}S_i(z)=S(z)$ exists.
  \end{lem}
\begin{proof}
Let $k$ be fixed and consider $U_k=(G_k\circ \dots\circ
G_1)^{-1}(\frac1{2C_2} D(r_2))$, where $C_2$, $ r_2$ are as in
Lemma \ref{correctionmaps}(c). Let $\Psi_k$ as in Lemma
\ref{correctionmaps}(c). Since the convergence in the lemma is
uniform in $D(r_2)$,   and $G^{-1}_{k+1}\circ
G^{-1}_{k+2}\circ\cdots\circ G^{-1}_{k+l}\circ T_{k+l}\circ
F_{k+l}\circ\cdots\circ F_{k+2}\circ F_{k+1}$ and $\Psi_k$ are
injective, we can find $0<\rho<r_2$ and $l_0$ such that
\begin{equation}\label{1}
G^{-1}_{k+1}\circ G^{-1}_{k+2}\circ\cdots\circ G^{-1}_{k+l}\circ
T_{k+l}\circ F_{k+l}\circ\cdots\circ F_{k+2}\circ
F_{k+1}(D(\rho))\supset  D(\frac1{2C_2}r_2)
\end{equation}
if $l\ge l_0$. Hence for every  $l\ge l_0$ we have: for every $z\in U_k$, there
exists $\zeta_l\in D(\rho)$ such that
\begin{equation}\label{2}
G_1^{-1}\circ \cdots\circ  G_k^{-1}\circ  G^{-1}_{k+1}\circ
\cdots\circ  G^{-1}_{k+l}\circ  T_{k+l}\circ  F_{k+l}\circ  \cdots
\circ F_{k+2}\circ  F_{k+1}(\zeta_l)=z.
\end{equation}

Hence
\begin{equation}\label{3}
\begin{split}
S_{k+l}(z)&=\Phi_{k+l}\circ T_{k+l}^{-1}\circ G_{k+l}\circ \cdots \circ G_1(z)\\
&=\Phi_{k+l}\circ F_{k+l}\circ  \cdots\circ  F_{k+2}\circ  F_{k+1}(\zeta_l)\\
&=\Phi_k(\zeta_l).
 \end{split}
\end{equation}

Take two such subsequences $\zeta_{l_j}$ and $\zeta_{l_j'}$ such that $\zeta_{l_j}\to
w$ and $\zeta_{l_j'}=w'$ as $j\to \infty$ with $w, w'\in \ol{D(\rho)}$. Since the
convergence in Lemma \ref{correctionmaps} is uniform, by (\ref{2}) we have
$$
G_1^{-1}\circ \cdots\circ  G_k^{-1}\circ \Psi_k(w)=z=G_1^{-1}\circ
\cdots\circ  G_k^{-1}\circ \Psi_k(w').
$$
Hence we must have $w=w'$ and so $\zeta_l\to w$ as $l\to\infty$. By (\ref{3}), we
have
\begin{equation}\label{4}
\begin{split}
\lim_{l\to\infty}S_{k+l}(z)&=\Phi_k(w)\\
&=\Phi_k\circ \Psi_k^{-1}\circ G_k\circ \cdots\circ  G_1(z).
\end{split}
\end{equation}

Hence $S=\lim_{i\to\infty} S_i$ exists on $U_k$. By Lemma \ref{correctionmaps},
$\bigcup_k U_k=\cn$, from this the lemma follows.
\end{proof}
\begin{lem}\label{nondegenerate}
 $S$ is a local biholomorphic map from $\cn$ into $M$.
\end{lem}
\begin{proof} This follows from (\ref{4}) immediately.
\end{proof}

\begin{lem}\label{exhaustion} For any $\e>0$,
$\bigcup_k\Phi_k(D(\e))=M$.
\end{lem}
\begin{proof}   Since the Ricci curvature
of $g(0)$ is positive, for any $R>0$, we have $R_\ijb(x,0)\ge ag_\ijb(x,0)$ for
some $a>0$ for all $x\in B_0(p,R)$ which is the geodesic ball with respect to
$g(0)$. Let $\lambda_i(x,t)$ be the eigenvalues of $R_\ijb(x,t)$ with respect to
$g(t)$. Then $\lambda_i(x,t)\le C$ for some constant $C$ independent of $x$ and
$t$. On the other hand, by Proposition \ref{s2p2}, we have
$\prod_{i=1}^n\lambda_i(x,t)\ge a^n$ for $x\in B_0(p,R)$. Hence there exists
$b>0$ independent of $t$ such that
$$
R_\ijb(x,t)\ge bg_\ijb(x,t)
$$
for all $t\ge 0$ and $x\in B_0(p,R)$. By the K\"ahler-Ricci flow equation, we
conclude that
$$
B_0(p,R)\subset B_t(p,e^{-\frac {b+\kappa}2t}R).
$$
From this and   Lemma \ref{s3l1} the lemma follows.
\end{proof}

\begin{lem}\label{surjective} $S$ is surjective.
\end{lem}
\begin{proof}
 From the proof of Lemma \ref{limit}, we conclude that
$$
S(\cn)\supset \Phi_k\circ \Psi_k^{-1}(D(\frac1{2C_2}r_2))
$$
for all $k$.  From this, Lemma \ref{exhaustion} and the proof of Lemma
 \ref{correctionmaps} (c), it is easy to
see that $S(\cn)=M$.
\end{proof}

\begin{lem}\label{injective} $S$ is injective.
\end{lem}
\begin{proof}
 Suppose there exists distinct  $z_1, z_2\in \cn$
such that $S(z_1)=S(z_2)=q$. Let $\sigma(\tau)$, $0\le \tau\le1$
be the line segment from $z_1$ to $z_2$. Let $\gamma=S\circ
\sigma$. Then $\gamma$ is a smooth closed curve in $M$ starting
from $q$. Since $M$ is simply connected, we can find a smooth
homotopy $\alpha(s,\tau)$,  $0\le s, \tau\le1$, with
$\alpha(0,\tau)=\gamma(\tau)$, $\alpha(1,\tau)=q$ (the constant
map), and $\alpha(s,0)=\alpha(s,1)=q$ for all $s$.

By Lemma \ref{correctionmaps}, there exists $0<\rho<r_2$ and
$\eta>0$ which are independent of $i$ such that $\Psi_i^{-1}$ is
defined on $D(\eta)$ and
\begin{equation}\label{injective1}
   \Psi_i^{-1}(D(\eta))\subset  D(\frac18\rho)
\end{equation}
For  $\rho>0$ let  $\rho_1>0$ be such that Lemma \ref{s3l2} is
true.   By the proof of Lemma \ref{exhaustion}, there exists $k_0$
such that if $k\ge k_0$, then $q\in B_k(p,\rho_1)$ and
$L_k(\alpha(s,\cdot))< \rho_1$ for all $s$, where $B_k(p,\rho_1)$
is the geodesic and $L_k$ is the length with respect to $g(kT)$.
 We may also assume that $\sigma \subset G_1^{-1}\circ \cdots\circ
 G_k^{-1}(D(\eta))$ provided $k_0$ is large enough.

Now fix $k\ge k_0$. Let $\wt\gamma=\Psi_k^{-1}\circ G_k\circ
\cdots\circ  G_1\circ \sigma$. Then $\wt\gamma\subset
D(\frac\rho8)$ by (\ref {injective1}) and it is a lift of $\gamma$
by $\Phi_k$ by (\ref{4}). Moreover, if $\wt\gamma(0)=w_1$, and
$\wt\gamma(1)=w_2$, then $w_1\neq w_2$ because $G_i$ and $\Psi_k$
are injective. Since $\Phi_k(w_1)=\Phi_k(w_2)=q\in B_k(p,\rho_1)$,
by Lemma \ref{s3l2}, for any $s$, there is a lift
$\wt\beta_s(\tau)$ of $\alpha(s,\tau)$ by $\Phi_k$ such that
$\wt\beta_s(0)=w_1$ and $\wt\beta_s\subset D(\frac\rho2)$.

We claim that $\wt\beta_s(1)=w_2$. Let $\e>0$ be such that
$\Phi_k$ is a biholomorphism onto its image when  restricted on
$B_e(w_2,\e)$. For such $\e>0$, let $\delta>0$ be as in Lemma
\ref{s3l3}. On the other hand, let $\xi>0$ be such that
\begin{equation}\label{uniformcont}
d_k(\alpha(s_1,\tau),\alpha(s_2,\tau))<\delta
\end{equation}
for all $\tau$, if  $|s_1-s_2|\le \xi$.

Hence by  Lemma \ref{s3l3},   if $0\le s\le \xi$ then
$\wt\beta_s(1)\in B_e(w_2,\e)$. But
$\Phi_k(\wt\beta_s(1))=\alpha(s,1)=q$, and $\Phi_k$ is injective
on $ B_e(w_2,\e)$.  Thus we have $\wt\beta_s(1)=w_2 $ for $0\le
s\le \xi$. In particular, $\wt\beta_\xi(1)=w_2$. By
(\ref{uniformcont}) and Lemma \ref{s3l3}, we can argue as before
and conclude that $\wt\beta_{2\xi}(1)=w_2$. Continue in this way,
we have that $\wt\beta_1(1)=w_2$. On the other hand, $\wt\beta_1$
is a lift of $\alpha(1,\cdot)$. Hence $\Phi_k(\wt\beta_1(\tau))=q$
for all $\tau$. Since $w_1\neq w_2$, there is $\tau$ with
$|\wt\beta_1(\tau)-w_2|=\frac\e2$. This is impossible because
$\Phi_k$ is a injective on $B_e(w_2,\e)$.
\end{proof}
Now Theorem  \ref{s1t3} follows from Lemmas \ref{nondegenerate},
\ref{surjective}, and \ref{injective}.

\bibliographystyle{amsplain}

\end{document}